\documentclass[11pt]{amsart}
\usepackage{amssymb}

\usepackage[T1]{fontenc} 
\usepackage[utf8]{inputenc} 

\usepackage[unicode,colorlinks,
	citecolor=blue,
	urlcolor=blue,
	linkcolor=blue
]{hyperref} 

\usepackage{soul}

\usepackage{tikz}

\theoremstyle{plain}
\newtheorem{thm}{Theorem}[section]

\newtheorem{lem}[thm]{Lemma}

\theoremstyle{definition}
\newtheorem{defn}[thm]{Definition}

\theoremstyle{remark}
\newtheorem{remark}[thm]{Remark}

\numberwithin{equation}{section}

\renewcommand{\epsilon}{\varepsilon}
\renewcommand{\phi}{\varphi}

\renewcommand{\preceq}{\preccurlyeq}

\newcommand{\lb}{\lbrace}
\newcommand{\lk}{\lbrack}
\newcommand{\rb}{\rbrace}
\newcommand{\rk}{\rbrack}

\newcommand{\ra}{\rightarrow}

\begin{document}

\title[Introenumerability, Autoreducibility, and Randomness]{Introenumerability, Autoreducibility,\\ and Randomness}

\author{Ang Li}

\date{\today}

\begin{abstract}
 We define $\Psi$-autoreducible sets given an autoreduction procedure $\Psi$. Then, we show that for any $\Psi$, a measurable class of $\Psi$-autoreducible sets has measure zero. Using this, we show that classes of cototal, uniformly introenumerable, introenumerable, and hyper-cototal enumeration degrees all have measure zero.  
 
 By analyzing the arithmetical complexity of the classes of cototal sets and cototal enumeration degrees, we show that weakly 2-random sets cannot be cototal and weakly 3-random sets cannot be of cototal enumeration degree. Then, we see that this result is optimal by showing that there exists a 1-random cototal set and a 2-random set of cototal enumeration degree. For uniformly introenumerable degrees and introenumerable degrees, we utilize $\Psi$-autoreducibility again to show the optimal result that no weakly 3-random sets can have introenumerable enumeration degree. We also show that no 1-random set can be introenumerable.
\end{abstract}

\maketitle

\section{Introduction}

In 1959, Friedberg and Rogers \cite{enumerationreducibility} introduced enumeration reducibility. A set $A\subseteq\omega$ is \emph{enumeration reducible} to another set $B\subseteq\omega$ if there is a c.e.\ set $W$ such that $A=\lb x:(\exists y)\langle x,y\rangle\in W\mathrm{\ and\ }D_y\subseteq B\rb$, where $\lb D_y\rb_{y\in\omega}$ gives a computable listing of all finite sets. We call the c.e.\ set $W$ that witnesses this reduction an enumeration operator and write $A=W(B)$. The degree structure induced by enumeration reduction $\leq_e$ consists of the enumeration degrees. We can identify subsets of $\omega$ with infinite strings in the Cantor space $2^{\omega}$. Therefore, we can consider the measure of different classes of enumeration degrees (often abbreviated by e-degrees), including cototal e-degrees, uniformly introenumerable e-degrees, introenumerable e-degrees, and hyper-cototal e-degrees. 
\par
Given a set $A$ of natural numbers and any number $n$, we may ask whether the membership of $n$ in $A$ can be determined using the oracle $A$ without asking ``is $n$ in $A$''. If so, $A$ has a kind of self-reducibility. The notion of autoreducibility introduced by Trakhtenbrot \cite{autoreducibility} in 1970 is a formalization of this idea. A set $A$ is said to be \emph{autoreducible} if there is a Turing functional $\Phi$ such that for any $n$, $A(n)=\Phi^{A-\lb n\rb}(n)$. We will generalize the autoreduction notion by defining $\Psi$-autoreducibility for any autoreduction procedure $\Psi$, which is a function from $\omega\times2^{\omega}$ to $\lb0,1\rb$. The classes of enumeration degrees mentioned above all have natural autoreducibility by replacing the Turing functional with different autoreduction procedures. Next, we will show that any measurable class of $\Psi$-autoreducible sets has measure zero for any $\Psi$. Then, we use this property of classes of $\Psi$-autoreducible sets to show that the classes of above e-degrees all have measure zero.
\par
Intuitively, given a set $A\subseteq\omega$, or equivalently an infinite string in $2^{\omega}$, it is random if it is hard to compress or one cannot predict the next bit or it has no rare properties. In 1966, Martin-L\"of introduced a randomness notion using the latter idea that a random set is in no effective measure zero set in \cite{random}. An infinite string $A\in2^{\omega}$ is \emph{Martin-L\"of random} (\emph{n-random}), if $A$ is not in $\bigcap_m G_m$ where $\lb G_m\rb_{m\in\omega}$ is any uniformly $\Sigma^0_1$ ($\Sigma^0_n$ respectively) sequence of open sets such that the measure of each $G_m$ is smaller than $2^{-m}$. A set $A$ is \emph{weakly $n$-random} if $A$ avoids all $\Pi^0_n$ classes. Generally, a set is random if it avoids a particular kind of null classes. Such null classes can be arithmetical as above or even go beyond arithmetical. 
\par
Since our classes of e-degrees have measure zero, sufficiently random sets must avoid such measure zero classes. Therefore, we can ask questions about what level of randomness the above sets or e-degrees can reach, and what level of randomness the above sets or e-degrees must avoid. We answer such questions for cototal sets, cototal e-degrees, uniformly introenumerable sets, uniformly introenumerable e-degrees, introenumerable sets, and introenumerable e-degrees. For references for randomness notions, see \cite{downey} or \cite{Nies}.
\par
We start by giving the definitions of the sets and e-degrees we mentioned. First, a set $A$ is \emph{total} if $\overline{A}\leq_e A$. It is named total because the degree of a total set is the degree of the graph of a total function. In \cite{Cototal}, the notion of cototality is given by reversing the relationship between $A$ and $\overline{A}$.

\begin{defn}
A set $A$ is \emph{cototal} when $A\leq_e\overline{A}$.
\end{defn}

\begin{defn}
An infinite set $X$ is \emph{uniformly introenumerable} if there is an enumeration operator $\Gamma$ such that for every infinite subset $Y$ of $X$, $\Gamma(Y)=X$.
\end{defn}

In \cite{uie}, Jockush introduced the notion of uniform introenumerability. The definition of uniform introenumerability we give here is slightly different by using an enumeration operator instead of a c.e.\ operator, though the two definitions were shown to be equivalent in \cite{Introreducible} by Greenberg et al. Recently, Goh et al. \cite{GJ-GMS} also showed that Jockush's notion of (non-uniform) introenumerability is equivalent to the following notion:

\begin{defn}
An infinite set $X$ is \emph{introenumerable} if, for every infinite subset $Y$ of $X$, there is an enumeration operator $\Gamma$ such that $\Gamma(Y)=X$.
\end{defn}

In \cite{Sanchis}, Sanchis introduced a reduction that is related to hyperarithmetical reduction and only uses positive information about membership in the set: 

\begin{defn}\label{Def_hypc}
Let $A$ and $B$ be sets such that, for some c.e.\ set $W$, the following
relation holds:
$x\in B$ if and only if $$(\forall f\in\omega^{<\omega})(\exists n,y)\lk\langle f\upharpoonright n,x,y\rangle\in W\wedge D_y\subseteq A\rk.$$
Then we say that $B$ is \emph{hyper-enumeration} reducible to $A$ and write this
relation: $B\leq_{he}A$.
\end{defn}

\begin{defn}
A is called \emph{hyper-cototal} if $A \leq_{he}\overline{A}$. 
\end{defn}

\begin{thm}\label{Thm_rel}
The relationship of enumeration degrees of the above notions is the following:
\begin{align*}
    \textbf{Cototal} & \rightarrow\textbf{Uniformly Introenumerable}\\ & \dashrightarrow\textbf{Introenumerable}\ra\textbf{Hyper-cototal}.
\end{align*}

\end{thm}

\begin{remark}
The solid arrows are strict. For proof of the first arrow, see \cite{McCarthy}. The third arrow and the strictness of the first arrow are proved in \cite{GJ-GMS} by Goh et al. It is still unknown whether there is a set of introenumerable e-degree that does not have uniformly introenumerable e-degree.
\end{remark}

\section{Measure of Classes with Autoreduction}

In this section, we define $\Psi$-autoreducible sets given an autoreduction procedure $\Psi$ and show that any measurable class of $\Psi$-autoreduction sets has measure zero. Next, we apply the autoreducibility of hyper-cototal e-degrees to show that the measure of the class of such e-degrees is zero.

\begin{defn}
Given a function $\Psi:\omega\times2^{\omega}\ra\lb0,1\rb$, A set $A$ is $\Psi$-\emph{autoreducible} if and only if $$(\forall n)\lk A(n)=\Psi(n,A-\lb n\rb)\rk.$$ Here, we say that the function $\Psi$ is an autoreduction procedure.
\end{defn}

Next, to show that the measure of a class of $\Psi$-autoreducible sets is zero, we use the Lebesgue density theorem.
\begin{thm}\label{Thm_auto}
Fix an autoreduction procedure $\Psi$, a measurable class $S$ of $\Psi$-autoreducible sets has measure zero.
\end{thm}

\begin{proof}
Suppose a class $S$ of $\Psi$-autoreducible sets has positive measure. By the Lebesgue density theorem, for any $\epsilon>0$, there is a string $\sigma\in 2^{<\omega}$ such that $\frac{\mu(S\cap \lk\sigma\rk)}{\mu(S)}\geq1-\epsilon$. Fix $\epsilon=\frac{1}{4}$ along with the corresponding string $\sigma$. Consider an $n\in\omega$ larger than $|\sigma|$. Define subsets $P_i$ $(i=0,1)$ of $S$ as follows: $$P_i=\lb X\in S:\Psi(n,X-\lb n\rb)=i\rb.$$ Since $P_0$ and $P_1$ partition $S$, one of them must have the following relative measure: $\frac{\mu(P_i\cap\lk\sigma\rk)}{\mu(S)}\geq\frac{1-\epsilon}{2}=\frac{3}{8}$. Without loss of generality, assume that such subset is $P_0$. Now, consider the set $$P_2=\lb\hat{X}:X\in P_0, \hat{X}(n)=1, (\forall i\not=n)\lk X(i)=\hat{X}(i)\rk\rb.$$ Notice that if $x\in P_0$, $X(n)=0$. So, $P_2$ also has relative measure $\frac{\mu(P_2\cap\lk\sigma\rk)}{\mu(S)}\geq\frac{3}{8}>\frac{1}{4}$. Therefore, $\frac{\mu(P_2\cap S\cap \lk\sigma\rk)}{\mu(S)}>0$. So, $P_2\cap S$ is not empty. For any $Y\in P_2\cap S$, $\Psi(n,Y-\lb n\rb)=0\not=1=Y(n)$. This is a contradiction. Therefore, $S$ has measure zero.
\end{proof}

\begin{remark}
In this theorem, the assumption that the class $S$ is measurable is necessary. Consider the finite difference equivalence classes: two sets $A$ and $B$ are in the same equivalence class if and only $(A-B)\cup(B-A)$ is finite. Now, we can define a class $S_0$ that contains exactly one element from each of the equivalence classes. It is not difficult to see that $S_0$ is not measurable. We can define a function $\Psi_0$ such that if $A\in S_0$ and $n\in\omega$, then $\Psi_0(n,A-\lb n\rb)=A(n)$. It is well-defined because, for any $B\in S_0$ and $B-\lb n\rb=A-\lb n\rb$, $\Psi_0(n,B-\lb n\rb)$ has to equal $A(n)$ by the definition of $S_0$. Therefore, $S_0$ is a class consisting of $\Psi_0$-autoreducible sets that does not have measure zero since it is not measurable.
\end{remark}

Now we use the above theorem to show that the measure of the class of hyper-cototal e-degrees is zero. First, we discuss the autoreducibility of hyper-cototal sets. 

\begin{lem}\label{Lem_hypauto}
For every hyper-cototal set $A$, there is a $\Psi$ such that $A$ is $\Psi$-autoreducible.
\end{lem}

\begin{proof}
Suppose $A$ is hyper-cototal and there is some hyper-enumeration operator $\Delta$ such that $A=\Delta(\overline{A})$. When $n\in A$, $\overline{A}\subseteq\overline{A-\lb n\rb}$. Therefore, $n\in\Delta^{\overline{A}}\subseteq\Delta^{\overline{A-\lb n\rb}}$. When $n\not\in A$, $n\not\in\Delta^{\overline{A}}=\Delta^{\overline{A-\lb n\rb}}$. So, $A(n)=\Delta^{\overline{A-\lb n\rb}}(n)$. Then, we can define $\Psi(n,X):=\Delta^{\overline{X}}(n)$.
\end{proof}

In fact, each set of hyper-cototal degree is $\Psi$-autoreducible for some autoreduction procedure $\Psi$ as well. 

\begin{lem}
Any set in the class of hyper-cototal e-degrees is a hyper-cototal set.
\end{lem}

\begin{proof}
In \cite{Sanchis}, Sanchis proved that If $A\leq_e B$, then $A\leq_{he} B$ and $\overline{A}\leq_{he}\overline{B}$. Suppose $A$ has hyper-cototal e-degree and $A\equiv_e B$, where $B$ is a hyper-cototal set. Then, $A\equiv_{he} B\leq_{he}\overline{B}\equiv_{he}\overline{A}$.
\end{proof}

Next, in order to apply Theorem \ref{Thm_auto} to show that the measure of the classes of hyper-cototal e-degrees is 0, we first need to show that the class of hyper-cototal e-degrees is measurable by analyzing the arithmetical complexity of $$\lb A:A\leq_{he}\overline{A}\rb=\bigcup_{\Gamma}\lb A:(\forall n)\lk n\in A\ra n\in \Gamma^{\overline{A}}\wedge n\not\in A\ra n\not\in\Gamma^{\overline{A}}\rk\rb.$$ Notice that $n\in\Gamma^{\overline{A}}$ and $n\not\in\Gamma^{\overline{A}}$ are $\Pi_1^1$ and $\Sigma_1^1$ respectively for a hyper-enumeration operator $\Gamma$ by Definition \ref{Def_hypc}. So, the class of hyper-cototal e-degrees is the difference of two $\Pi_1^1$ classes. Recall that $\Pi^1_1$ sets are measurable. Therefore, the class of hyper-cototal e-degrees is measurable. Now, we use the results from above to see that the class of hyper-cototal e-degrees has measure zero.

\begin{lem}\label{Cor_mea}
The classes of hyper-cototal, introenumerable, uniformly introenumerable, and cototal e-degrees all have measure zero.
\end{lem}

\begin{proof}
Suppose the class of hyper-cototal e-degrees has positive measure. Because there are only countably many hyper-enumeration operators, there exists a $\Gamma$ such that the class of hyper-cototal e-degrees witnessed by this operator has positive measure. However, any set in this class would be $\Gamma$-autoreducible by Lemma \ref{Lem_hypauto}. Now, applying Theorem \ref{Thm_auto} gives us a contradiction.
By the relationship between the e-degrees mentioned above in Theorem \ref{Thm_rel}, we see that the measure of these classes are all zero.
\end{proof}

\section{Bounds of Randomness}

Notice that, for any class of measure zero, sufficiently random sets avoid it. So, we now discuss what level of randomness these e-degrees could and could not have. In this section, all necessary background knowledge of randomness is from Nies' book \cite{Nies}. We first discuss the class of cototal sets and the class of cototal e-degrees.

\begin{thm}\label{Thm_W2R}
Weakly 2-random sets are not cototal.
\end{thm}

\begin{proof}
The class of cototal sets $\lb A:A\leq_e \overline{A}\rb$ is defined by 
    \begin{align*}
        \bigcup_e\lb A:A=\Gamma_e^{\overline{A}}\rb= & \bigcup_e\lb A:\forall n\lk  n\in A\ra(\exists D_y\subseteq \overline{A})\lk\langle n,y\rangle\in\Gamma_e\rangle\rk\\ 
        & \wedge n\not\in A\ra (\forall y)\lk\langle n,y\rangle\in\Gamma_e \ra D_y\cap A\not=\emptyset\rk\rk\rb,
    \end{align*}
where $\Gamma_e$'s are enumeration operators. Therefore, the class of cototal sets is a union of $\Pi^0_2$ classes. By Lemma \ref{Cor_mea}, all such classes have measure zero. Because any weakly 2-random set avoids all null $\Pi^0_2$ classes, weakly 2-random sets are not cototal.
\end{proof}

To see that weak 2-randomness is optimal, we show that the $1$-random Chaitin's $\Omega$ is a cototal set.

\begin{thm}
There exists a 1-random cototal set.
\end{thm}

\begin{proof}
Because $\Omega$ is left-c.e., there is a non-descending computable sequence $\lb q_n\rb$ of rationals such that $\Omega=\lim_{n\ra\infty}q_n$. For any enumeration of $\overline{\Omega}$, we can enumerate $\Omega$ using this computable sequence. First, to determine whether $0$ is in $\Omega$ or not, either for some $n$, we see the dyadic expansion of $q_n$ starts with $1$ or we see $1$ enter $\overline{\Omega}$. Only for the first case, we enumerate $0$ in $\Omega$. Then, we can iteratively do this process for each nature number in order. Eventually, we obtain an enumeration of $\Omega$. Therefore, $\Omega\leq_e\overline{\Omega}$.
\end{proof}

For the class of cototal e-degrees, we first discuss what level of randomness is enough to avoid them.

\begin{thm}
Weakly 3-random sets do not have cototal e-degree.
\end{thm}

\begin{proof}
Notice that the class of cototal e-degrees defined by an enumeration operator $\Gamma_e$ is
\begin{align*}
    \lb A:A=\Gamma_e^{\overline{K_A}}\rb= \lb A:(\forall n) \lk & n\in A\ra(\exists y)\lk\langle n,y\rangle\in\Gamma_e\ra D_y\cap K_A=\emptyset\rk\\ \wedge  & n\not\in A\ra (\forall y)\lk \langle n,y\rangle\in\Gamma_e\ra D_y\cap K_A\not=\emptyset\rk\rk\rb.
\end{align*}
 Since $D\cap K_A=\emptyset$ and $D\cap K_A\not=\emptyset$ are $\Pi^0_1$ and $\Sigma^0_1$ respectively, the class of cototal e-degrees defined by $\Gamma_e$ is $\Pi^0_3$. Since each of these classes is null, weakly 3-random sets avoid them all. So, we conclude that weakly 3-random sets do not have cototal e-degree.
\end{proof}

Next, we see that weak 3-randomness is optimal by showing that there is a 2-random set of cototal e-degree even though any cototal set cannot be weakly 2-random.

\begin{thm}\label{Thm_2R}
There exists a 2-random set of cototal e-degree.
\end{thm}

\begin{proof}
Consider Chaitin's $\Omega$ relativized to $\emptyset'$, i.e.\ $\Omega^{\emptyset'}$, which is 2-random. Let $L$ be $\lb q\in\mathbb{Q}_2: q<\Omega^{\emptyset'}\rb$. Then, $L\leq_e\Omega^{\emptyset'}\leq_e L\oplus\overline{L}$. Notice that $L$ is $\Sigma_2^0$. In \cite{Cototal}, it was shown that every $\Sigma^0_2$ set has cototal e-degree. So, there exists $M$ such that $M\equiv_e L$ and $\overline{M}\geq_e M$. Then, $\overline{\Omega^{\emptyset'}\oplus L\oplus M}\geq_e\overline{L}\oplus\overline{M}\geq_e\overline{L}\oplus M\equiv_e\overline{L}\oplus L\geq_e\Omega^{\emptyset'}\equiv_e\Omega^{\emptyset'}\oplus L\oplus L\equiv_e\Omega^{\emptyset'}\oplus L\oplus M$. Hence, we have a cototal set that is enumeration equivalent to $\Omega^{\emptyset'}$.
\end{proof}

In the proofs above, we did not use autoreducibility since it is enough to analyze the arithmetical complexities of the class of cototal sets and the class of cototal e-degrees to show the optimal level of randomness the sets in these classes must avoid. However, a similar analysis would not work for the classes of (uniform) introenumerable sets or e-degrees. We can verify the complexity of the collection of uniformly introenumerable e-degrees: 
\begin{align*}
    \bigcup_e\lb A:\exists i,m\forall B & \lk\forall a \lk a\in A\leftrightarrow\exists b\lk\langle a,b\rangle\in\Gamma_m\wedge D_b\subseteq \Gamma_i(A)\rk\rk\wedge \\
    & \lk B\subseteq\Gamma_i(A)\wedge\lk\forall p\in B\exists q>p \rk\ra \\
    \forall t & \lk t\in \Gamma_i(A)\leftrightarrow \exists s\lk \langle t,s\rangle\in\Gamma_e\wedge D_s\subseteq B\rk\rk\rk\rk\rb.
\end{align*}
 This is $\Pi^1_1$. We suspect that the class of uniformly introenumerable e-degrees is $\Pi^1_1$-complete. This was shown to be true for the class of uniformly introreducible sets in \cite{Introreducible}. Assuming that there is no simpler definition, the analysis we used for cototal e-degrees would not work. Instead, for each set $A$ of uniform introenumerable e-degree, we show $\Psi$-autoreducibility for some autoreduction procedure $\Psi$ so that we can apply Theorem \ref{Thm_auto} again.

\begin{thm}\label{Thm_W3R}
Weakly 3-random sets do not have uniformly introenumerable e-degree.
\end{thm}

\begin{proof}
We will show that uniformly introenumerable e-degrees are contained in a countable union of measure zero $\Pi^0_3$ classes.
To do this, we show that each set $A$ of uniformly introenumerable e-degree is $\Psi$-autoreducible for some $\Psi$. Since $A$ has uniformly introenumerable e-degree, there is a set $B$, enumeration operators $\Phi$, $\Gamma$, and $\Delta$ such that $A=\Delta(B),B=\Phi(A)$, and for any infinite subset $C$ of $B$, $\Gamma(C)=B$. Let $$\Psi(n,Z)=\begin{cases} 1 & n\in\Delta(\Gamma(\Phi(Z))) \\ & \text{ or } \Phi(Z)\text{ is finite,}\\ 0 & \text{otherwise}. \end{cases}$$ 
Note that $n$ has to be in $A$ when $\Phi(A-\lb n\rb)$ is finite. So, $A$ is $\Psi$-autoreducible. Now we consider the class of $\Psi$-autoreducible sets: 
\begin{align*}
    \lb D:\forall n\lk & \lk n\in D\rightarrow n\in\Delta(\Gamma(\Phi(D-\lb n\rb)))\\ & \vee(\exists p\forall t>p) \lk t\not\in\Phi(D-\lb n\rb)\rk\rk\\ & \wedge \lk n\not\in D\rightarrow(\forall q\exists s>q)\lk s\in\Phi(D-\lb n\rb)\rk\\ & \wedge n\not\in\Delta(\Gamma(\Phi(D-\lb n\rb)))\rk\rk\rb.
\end{align*}
 This is a $\Pi^0_3$ class. By Theorem \ref{Thm_auto}, this is a null class. Because weakly 3-random sets cannot be in any $\Pi^0_3$ null class, weakly 3-random sets do not have uniformly introenumerable e-degree.
\end{proof}

Meanwhile, there also exists $2$-random uniformly introenumerable e-degrees because of Theorem \ref{Thm_2R} and the fact that every set of cototal e-degree has uniform introenumerable e-degree.

\par
With more work, the previous result can be improved to show that weakly 3-random sets do not have introenumerable e-degree either. 

\begin{thm}\label{Thm_intro}
No weakly $3$-random set has introenumerable e-degree.
\end{thm}

\begin{proof}
Suppose a weakly 3-random set $A$ has introenumerable e-degree. Let $B$ be an introenumerable set such that there are enumeration operators $\Phi$ and $\Delta$ with $A=\Delta(B)$ and $B=\Phi(A)$. For a contradiction, we define $C=\bigcup_ic_i$ as an infinite subset of $B$ such that $\Gamma_i(C)\not=B$ for any enumeration operator $\Gamma_i$ (here we identified strings $c_i$ with corresponding sets). When we are constructing $C$, we also define a set $D_i$ at each stage $i$. Let $c_0=\emptyset$ and $D_0=\emptyset$. Suppose $c_i$ and $D_i$ have been defined. By inductive assumption, $\Phi(A-D_i)$ is infinite. First, we consider whether there is an extension $e$ of $c_i$ such that $e\preceq c_i\Phi(A-D_i)\upharpoonright\lk|c_i|,\infty)$, and $\Gamma_i(e)-B\not=\emptyset$. If so, we define $c_{i+1}$ to be the least such $e$ that contains at least one more element than $c_i$ and $D_{i+1}=D_i$. If not, we consider whether there is an extension $e$ of $c_i$ such that for some $n$, $\Phi(A-D_i\cup\lb n\rb)$ is infinite, $e\preceq c_i\Phi(A-D_i)\upharpoonright\lk|c_i|,\infty)$, and $\Gamma_i(e\Phi(A-D_i\cup\lb n\rb)\upharpoonright\lk|e|,\infty))\subsetneq B$. If so, we define $c_{i+1}$ to be the least such $e$ that contains at least one more element than $c_i$, and $D_{i+1}=D_i\cup\lb n\rb$. If not, we can define $$\Psi(n,Z)=\begin{cases} 1 & n\in\Delta(\Gamma_i(c_i\Phi(Z-D_i)\upharpoonright\lk|c_i|,\infty))) \\ & \text{ or } \Phi(Z-D_i)\text{ is finite,}\\ 0 & \text{otherwise.} \end{cases}$$ similar to the proof in Theorem \ref{Thm_W3R}. Notice that $A$ is $\Psi$-autoreducible and the class of $\Psi$-autoreducible sets is $\Pi^0_3$. This is impossible because $A$ is weakly 3-random. This is a contradiction. Therefore, at least one of the two cases we considered has to be true. In this way, we obtain an infinite $C=\bigcup_ic_i\subseteq B$. Now we show that $\Gamma_i(C)\not=B$ for any $i$. For any $i$, if the first case we considered is true, then $\Gamma_i(C)$ contains an element not in $B$. If the second case is true, $\Gamma_i(C)\subseteq\Gamma_i(c_{i+1}\Phi(A-D_{i+1})\upharpoonright\lk|c_{i+1}|,\infty))\subsetneq B$.
\end{proof}

Again, by Theorems \ref{Thm_rel} and \ref{Thm_2R}, we conclude that there exists $2$-random introenumerable e-degree while there is no weakly $3$-random introenumerable e-degree. Next, we consider the class of uniformly introenumerable sets. We use the proof ideas of Proposition 8 given by Figueira, Miller, and Nies in \cite{Indifferent} that showed no random is autoreducible.  

\begin{thm}\label{Thm_No1}
No 1-random set is uniformly introenumerable.
\end{thm}

\begin{proof}
We will apply Schnorr's theorem. To do so, we will show that the initial segment of any uniformly introenumerable set $A$ can be compressed beyond any fixed constant. 
\par
Let $\Gamma$ be the enumeration operator such that $\Gamma(B)=A$ for any infinite subset $B$ of $A$. For each $m$, there is a least $n_m$ such that $n_m>n_p$ for any $p<m$ and $\Gamma_{n_m}(0^mA\upharpoonright\lk m,n_m))\upharpoonright m=A\upharpoonright m$ since $A-\lb0,1,...,m-1\rb$ is an infinite subset of $A$. Let $c_m$ be the number of $1$'s in the string $A\upharpoonright m$. 
\par
Now we define a prefix-free machine $M$ that outputs $A\upharpoonright n_m$ with input $\gamma=0^{|\sigma|}1\sigma0^{|\tau|}1\tau A\upharpoonright\lk m,n_m)$, where $\sigma,\tau$ are binary strings corresponding to $m,c_m$. $M$ first obtains the length of $\sigma$ by reading until the first $1$ and then obtains the number $m$ by reading $|\sigma|$ many bits after the first $1$. Next, $M$ can find out $c_m$ in the same way by reading the input until $\tau$. Now, $M$'s read head keeps on moving forward to read $A\upharpoonright\lk m,n_m)$ bit by bit to do the enumeration of $\Gamma(0^mA\upharpoonright\lk m,n_m))\upharpoonright m$ step by step to enumerate $A(x)$ for $x$ between $0$ and $m-1$ until $c_m$ many of such $A(x)$ is determined to be $1$, which means the other bits on $A\upharpoonright m$ are zeros. $M$ can output $A\upharpoonright n_m$ by concatenation. Therefore, $K(A\upharpoonright n_m)\leq^+ n_m-m+4\log(m)$. By Schnorr's theorem, $A$ is not 1-random.
\end{proof}

For introenumerable sets, we combine the methods used in Theorems \ref{Thm_intro} and \ref{Thm_No1}.

\begin{thm}
No 1-random set is introenumerable.
\end{thm}

\begin{proof}
Suppose there is a $1$-random introenumerable set $A$. We prove the theorem by constructing an infinite subset $B=\bigcup_i b_i$ of $A$ such that $\Gamma_i(B)\not=A$ for any enumeration operator $\Gamma_i$ (here we identified the strings $b_i$ with its corresponding set).
\par
Let $b_0=\emptyset$. Suppose we have already defined $b_i$. There are two possible cases. One of the two cases must hold for it to be $1$-random.
\par
First, We consider whether there is an $n$ such that $\Gamma_i(b_iA\upharpoonright\lk|b_i|,n))$ contains an element that is not in $A$. If so, we let $b_{i+1}=b_iA\upharpoonright\lk|b_i|,n)$. In this case, we have a finite extension $b_{i+1}$ of $b_i$ such that $b_{i+1}$ is a subset of $A$, and for any infinite extension $B$ of $b_{i+1}$, $\Gamma_{i}(B)$ has an element not in $A$. 
\par 
Second, if there is no such $n$ in the first case, we consider whether there is an $m$ such that $\Gamma_i(b_i0^mA\upharpoonright\lk|b_i|+m,\infty))\subsetneq A$. If so, we let $b_{i+1}=b_i0^m$. In this case, we have a finite extension $b_{i+1}$ of $b_i$ such that applying $\Gamma_i$ to $A$'s subset $b_{i+1}A\upharpoonright\lk|b_{i+1}|,\infty)$ does not output $A$.
\par
If one of the cases holds for every $i$, we can show that for any $i$, $\Gamma_i(B)\not=A$, contradicting introenumerability. If the first case holds for $i$, then for any extension $B_0$ of $b_{i+1}$, $\Gamma_i(B_0)\not=A$. If the first case does not hold, notice that $B$ is a subset of $B_1=b_i0^mA\upharpoonright\lk|b_i|+m,\infty)$. Then, $\Gamma_i(B)\subseteq\Gamma_i(B_1)\subsetneq A$.
\par
If neither cases hold for some $i$, we show that $A$ is not $1$-random using a method similar to the one used in the proof of the above theorem. For each $m$, there is a least $n_m$ such that $n_m>n_p$ for any $p<m$ and $$\Gamma_{i,n_m}(b_i0^mA\upharpoonright\lk|b_i|+m,n_m))\upharpoonright |b_i|+m=A\upharpoonright |b_i|+m$$ because the failure of the second case guarantees that eventually numbers in $A\upharpoonright |b_i|+m$ will be enumerated and no other numbers would be enumerated by the failure of the first case. Let $c_m$ be the number of $1$'s in the string $A\upharpoonright\lk|b_i|,|b_i|+m)$. Now we define a prefix-free machine $M$ that outputs $A\upharpoonright n_m$ with input $\gamma=0^{|\sigma|}1\sigma0^{|\tau|}1\tau A\upharpoonright\lk |b_i|+m,n_m)$, where $\sigma,\tau$ are binary strings corresponding to $m,c_m$. $M$ obtains $m,c_m$ in the same way as the proof above by reading until $\tau$. Then, $M$ obtains the first $|b_i|$ bits of $A$ using $\Gamma_i$. Next, its read head keeps on moving forward to read $A\upharpoonright\lk |b_i|+m,n_m)$ bit by bit to do the enumeration of $\Gamma_i(b_i0^mA\upharpoonright\lk|b_i|+m,n_m))$ step by step to enumerate $A(x)$ for $x$ between $|b_i|$ and $|b_i|+m-1$ until $c_m$ many of such $A(x)$ is determined to be $1$ and output $A\upharpoonright n_m$ by concatenation. Therefore, $K(A\upharpoonright n_m)\leq^+ n_m-m+4\log(m)$. By Schnorr's theorem, $A$ is not 1-random.
\end{proof}

%
%

\bibliographystyle{plain}
\bibliography{references} 

\begin{thebibliography}{10}

\bibitem{Cototal}
Uri Andrews, Hristo Ganchev, Rutger Kuyper, Steffen Lempp, Joseph Miller, Alexandra Soskova, and Mariya Soskova.
\newblock On cototality and the skip operator in the enumeration degrees.
\newblock {\em Transactions of the American Mathematical Society}, 372(3):1631--1670, 2019.

\bibitem{downey}
Rodney~G Downey and Denis~R Hirschfeldt.
\newblock {\em Algorithmic randomness and complexity}.
\newblock Springer Science \& Business Media, 2010.

\bibitem{Indifferent}
Santiago Figueira, Joseph Miller, and Andr{\'e} Nies.
\newblock Indifferent sets.
\newblock {\em Journal of Logic and Computation}, 19(2):425--443, 2009.

\bibitem{enumerationreducibility}
Richard~M. Friedberg and Hartley Rogers.
\newblock Reducibility and completeness for sets of integers.
\newblock {\em Mathematical Logic Quarterly}, 5:117--125, 1959.

\bibitem{GJ-GMS}
Jun~Le Goh, Josiah Jacobsen-Grocott, Joseph Miller, and Mariya Soskova.
\newblock Enumeration pointed trees.
\newblock In preparation.

\bibitem{Introreducible}
Noam Greenberg, Matthew Harrison-Trainor, Ludovic Patey, and Dan Turetsky.
\newblock Computing sets from all infinite subsets.
\newblock {\em Transactions of the American Mathematical Society}, 374(11):8131--8160, 2021.

\bibitem{uie}
Carl~G Jockusch.
\newblock Uniformly introreducible sets.
\newblock {\em The Journal of Symbolic Logic}, 33(4):521--536, 1969.

\bibitem{random}
Per Martin-L{\"o}f.
\newblock The definition of random sequences.
\newblock {\em Information and Control}, 9(6):602--619, 1966.

\bibitem{McCarthy}
Ethan McCarthy.
\newblock Cototal enumeration degrees and their applications to effective mathematics.
\newblock {\em Proceedings of the American Mathematical Society}, 146(8):3541--3552, 2018.

\bibitem{Nies}
Andr{\'e} Nies.
\newblock {\em Computability and randomness}, volume~51.
\newblock OUP Oxford, 2012.

\bibitem{Sanchis}
Luis~E Sanchis.
\newblock Hyperenumeration reducibility.
\newblock {\em Notre Dame Journal of Formal Logic}, 19(3):405--415, 1978.

\bibitem{autoreducibility}
Boris~A. Trakhtenbrot.
\newblock Autoreducibility.
\newblock {\em Dokl. Akad. Nauk SSSR}, 192(6):1224--1227, 1970.

\end{thebibliography}

\end{document}